\documentstyle[amscd]{amsart}

\theoremstyle{plain} \numberwithin{equation}{section}
\newtheorem{thm}{Theorem}[section]
\newtheorem{cor}[thm]{Corollary}
\newtheorem{lem}{Lemma}[section]
\newtheorem{prop}[thm]{Proposition}

\theoremstyle{remark}
\newtheorem*{rem}{Remark}

\begin{document}
\title[$({\Bbb Z}_2)^k$-actions with $w(F)=1$]
{\large \bf $({\Bbb Z}_2)^k$-actions with $w(F)=1$\\
\vskip .2cm {\tiny Dedicated to Professor Zhende Wu on his
seventieth birthday}}
 \author[Zhi L\"u]
{Zhi L\"u}
 \thanks{ Supported by  grants from NSFC (No. 10371020) and JSPS (No. P02299)}
\address{Institute of Mathematics, Fudan University, Shanghai,
200433, People's Republic of China.}
 \email{zlu@@fudan.edu.cn}
 \date{}
\begin{abstract}
Suppose that $(\Phi, M^n)$ is a smooth $({\Bbb Z}_2)^k$-action on
a closed smooth $n$-dimensional manifold such that all
Stiefel-Whitney classes of the tangent bundle to each connected
component of the fixed point set $F$ vanish in positive dimension.
This paper shows that if  $\dim M^n>2^k\dim F$ and each
$p$-dimensional part $F^p$ possesses the linear independence
property, then $(\Phi, M^n)$ bounds equivariantly, and in
particular, $2^k\dim F$ is the best possible upper bound of $\dim
M^n$ if $(\Phi, M^n)$ is nonbounding.
\end{abstract}

\maketitle {\large

 \section{Introduction}

Let $k$ be a positive integer. Suppose that $\Phi:({\Bbb
Z}_2)^k\times M^n\longrightarrow M^n$ is a smooth $({\Bbb
Z}_2)^k$-action on a closed smooth $n$-dimensional manifold. The
fixed point set $F$ of $(\Phi, M^n)$ consists of a union of closed
 submanifolds of different dimensions. By $\dim F$ we mean the dimension of
 the component  of $F$ of largest dimension.

In this paper, we are mainly concerned with the case in which all
Stiefel-Whitney classes of the tangent bundle to each connected
component of the fixed point set $F$ vanish in positive dimension,
i.e., $w(F)=1$, where $w$ denotes the total Stiefel-Whitney class.
For the case $k=1$, it was proved in \cite{l1} that if $\dim
M^n>2\dim F$, then $(\Phi, M^n)$ with $w(F)=1$ bounds
equivariantly. (Related  results with $k=1$ can  be found in
\cite{cf}, \cite{c}, and \cite{ks1}). Any involution is
equivariantly cobordant to an involution with the  property that
the $p$-dimensional part of the fixed set is connected. However,
for the case $k>1$, different components of the $p$-dimensional
part of the fixed set may have different normal representations.
This is just the key difficulty for the case $k>1$. In \cite{l2},
a linear independence condition for the fixed point set was
introduced. With the help of the condition, the argument can be
carried out without  the connectedness restriction for the fixed
point set, so that we may obtain the result in the general case.
Our main result is stated as follows.

\begin{thm}
Suppose that $(\Phi, M^n)$ is a smooth $({\Bbb Z}_2)^k$-action on
a closed smooth $n$-dimensional manifold such that each part $F^p$
of the fixed point set $F$ possesses the linear independence
property, and $w(F)=1$. If $\dim M^n>2^k\dim F$, then $(\Phi,
M^n)$ bounds equivariantly.
\end{thm}

\begin{rem} (1) When $k=1$, as shown in \cite{l1}, $2\dim F$ is the best
possible upper bound of $\dim M$ if the involution $(\Phi, M)$
with $w(F)=1$ doesn't bound. For the general case, Example 1 in
Section 3 will show that $2^k\dim F$ is still the best possible
upper bound of $\dim M$ if $(\Phi, M)$ doesn't bound.

(2) In Theorem 1.1, the condition that each part $F^p$ of the
fixed point set $F$ possesses the linear independence property is
necessary. This can be seen from Example 2 in Section 3.
\end{rem}

The method used here is the  formula given by Kosniowski and Stong
\cite{ks2}, which we will review in Section 2. The proof of
Theorem 1.1 will be finished in Section 3. The case $\dim
M=2^k\dim F$ will be discussed in Section 4.

The author expresses his gratitude to Professor R.E. Stong for his
valuable suggestions and providing the Example 1 in Section 3, and
also to Professor M. Masuda for helpful conversations.

\section{Kosniowski-Stong formula}

Let $G=({\Bbb Z}_2)^k$, and let $\text{Hom}(G, {\Bbb Z}_2)$ be the
set of all homomorphisms $\rho: G\longrightarrow {\Bbb Z}_2=\{\pm
1\}$, which consists of $2^k$ distinct homomorphisms. One agrees
to let $\rho_0$ denote the trivial element in $\text{Hom}(G,{\Bbb
Z}_2)$, i.e., $\rho_0(g)=1$ for all $g\in G$. Let
$EG\longrightarrow BG$ be the universal principal $G$-bundle,
where  $BG=EG/G=({\Bbb R}P^\infty)^k$ is the classifying space of
$G$. It is well-known that
$$H^*(BG;{\Bbb Z}_2)={\Bbb Z}_2[a_1,...,a_k]$$
with the $a_i$ one-dimensional generators. In particular, all
nonzero elements of $H^1(BG;{\Bbb Z}_2)\cong ({\Bbb Z}_2)^k $
consist of $2^k-1$ polynomials of degree one in ${\Bbb
Z}_2[a_1,...,a_k]$, i.e.,
$$a_1,..., a_k,$$
$$a_1+a_2,..., a_1+a_k, a_2+a_3,...,a_2+a_k,..., $$
$$\cdots$$
$$a_1+\cdots+a_{k-1},...,a_2+a_3+\cdots+a_k,$$
$$a_1+a_2+\cdots+a_k.$$
These polynomials of degree one correspond to all nontrivial
elements of $\text{Hom}(G,{\Bbb Z}_2)$ (note that actually
$H^1(BG;{\Bbb Z}_2)\cong \text{Hom}(G, {\Bbb Z}_2)$), and so for a
convenience,  they are denoted by $\alpha_\rho$ for $\rho\in
\text{Hom}(G,{\Bbb Z}_2)$ with $\rho\not=\rho_0$.  Also, one
agrees to let $\alpha_{\rho_0}=0$, the zero element of
$H^1(BG;{\Bbb Z}_2)$.

Let $X$ be a $G$-space. Then $X_G:=EG\times_GX$-----the orbit
space of the diagonal action on the product $EG\times X$-----is
the total space of the bundle $X\longrightarrow X_G\longrightarrow
BG$ associated to the universal principal bundle $G\longrightarrow
EG\longrightarrow BG$. The space $X_G=EG\times_GX$ is called the
{\em Borel construction} on the $G$-space $X$. Applying cohomology
with coefficients ${\Bbb Z}_2$ to the Borel construction $X_G$ on
the $G$-space $X$ gives the equivariant cohomology
$$H^*_G(X;{\Bbb Z}_2):=H^*(X_G;{\Bbb Z}_2).$$

Now let $(\Phi, M^n)$ be a smooth $G$-action on a smooth closed
manifold with nonempty fixed point set $F$, and let $\eta^{n_i}_i,
i=1,...,s$, be vector bundles with $G$-actions covering the action
$\Phi$ on $M$.  Then we have equivariant cohomologies
$H^*_G(M;{\Bbb Z}_2)$ and $H^*_G(F;{\Bbb Z}_2)$. They are all
$H^*(BG;{\Bbb Z}_2)$-modules; in particular, $H^*_G(F;{\Bbb Z}_2)$
is a free $H^*(BG;{\Bbb Z}_2)$-module (see, for example,
\cite{ap}).

Let
$$f(\alpha_\rho; x_1,...,x_n; x_1^1,...,x_{n_1}^1; \cdots;
x_1^s,...,x_{n_s}^s)$$ be a polynomial over ${\Bbb Z}_2$ which is
symmetric in each of the sets of variables $\{x_1,...,x_n\}$ and
$\{x_1^i,...,x^i_{n_i}\}, i=1,...,s$. In this polynomial, if we
let the $j$-th Stiefel-Whitney class of $M$ replace the $j$-th
elementary symmetric function in $x_1,...,x_n, \sigma_j(x)$, and
the $j_i$-th Stiefel-Whitney class of $\eta_i^{n_i}$ replace the
$j_i$-th elementary symmetric function in
$\{x_1^i,...,x_{n_i}^i\}$, then we obtain   a class in
$H^*(M;{\Bbb Z}_2)\otimes_{{\Bbb Z}_2}H^*(BG; {\Bbb Z}_2)$, which
may be evaluated on the fundamental homology class of $M$,  giving
an element
$$f(\alpha_\rho; x_1,...,x_n; x_1^1,...,x_{n_1}^1; \cdots;
x_1^s,...,x_{n_s}^s)[M]$$ in $H^*(BG;{\Bbb Z}_2)={\Bbb
Z}_2[a_1,...,a_k]$.  On the other hand, let $C$ be a connected
component of $F$ with $\dim C=p$. The irreducible real
$G$-representations are all one-dimensional and correspond to all
elements in Hom$(G, {\Bbb Z}_2)$. Actually, every irreducible real
$G$-representation has the form $\lambda_\rho: G\times{\Bbb
R}\longrightarrow {\Bbb R}$ by $\lambda_\rho(g, x)=\rho(g)x$.
$\lambda_{\rho_0}$ is the trivial representation determined by the
trivial element $\rho_0$. Then,  the normal bundle $\nu_C$ to $C$
in $M$ decomposes under the action of $G$ into the Whitney sum
$\bigoplus_{\rho\not=\rho_0}\nu_{C, \rho}$ of the subbundles
$\nu_{C,\rho}$, where each $\nu_{C, \rho}$ is the summand of
$\nu_C$ on which $G$ acts in the fibers via $\lambda_\rho$, and
let $\dim \nu_{C,\rho}=q_{C,\rho}$. Here we call  $\{q_{C,\rho}\}$
the {\em normal dimensional sequence of} $C$. Similarly, each
vector bundle $\eta_i^{n_i}$ restricted to $C$ decomposes under
the action of $G$ into the Whitney sum
$\bigoplus_{\rho}\eta_{\rho}^{n_i}$ of the subbundles
$\eta_{\rho}^{n_i}$, where $\eta_{\rho}^{n_i}$ is the subbundle on
which $G$ acts via $\rho$, and let $q_{\rho}^i=\dim
\eta_{\rho}^{n_i}$ so that $n_i=\sum_{\rho}q_{\rho}^i$. Now, in
the polynomial $f(\alpha_\rho; x_1,...,x_n; x_1^1,...,x_{n_1}^1;
\cdots; x_1^s,...,x_{n_s}^s)$, we replace $x_1,...,x_n$ by
$z_1,...,z_p$ and for all $\rho\not=\rho_0$, variables
$\alpha_\rho+y_{\rho}^i, 1\leq i\leq q_{C, \rho}$, and also
replace $x_1^i,...,x_{n_i}^i$ by the collection, for all $\rho$,
of $\alpha_\rho+v_{\rho}^{i,j}, 1\leq j\leq q_{\rho}^i$. Next, if
we replace the $j$-th elementary symmetric function in
$$\{z_1,...,z_p\}\text{  by  } w_j(C)=w_j(\tau_C),$$
$$\{y_{\rho}^1,...,y_{\rho}^{q_{C,\rho}}\}\text{ by }
w_j(\nu_{C,\rho}),$$
$$\{v_{\rho}^{i,1},...,v_{\rho}^{i, q_{\rho}^i}\}\text{ by }
w_j(\eta_{\rho}^{n_i}),$$ respectively, then the expression
$${{f(\alpha_\rho; z_1,...,z_p,
\alpha_\rho+y_{\rho}^1,...,\alpha_\rho+y_{\rho}^{q_{C,\rho}};
\alpha_\rho+v_{\rho}^{i,1},...,\alpha_\rho+v_{\rho}^{i,
q_{\rho}^i})}\over
{\prod_{\rho\not=\rho_0}\prod_{i=1}^{q_{C,\rho}}(\alpha_\rho+y_{\rho}^i)}}$$
 is a class in the localization $S^{-1}H^*_G(C;{\Bbb Z}_2)$ of the equivariant cohomology of $C$
 (here $S$ is the subset of $H^*(BG;{\Bbb Z}_2)$
 generated multiplicatively by nonzero elements in $H^1(BG;{\Bbb Z}_2)$), which may be
 evaluated on the fundamental homology class of $C$. This gives an
 element
$${{f(\alpha_\rho; z_1,...,z_p,
\alpha_\rho+y_{\rho}^1,...,\alpha_\rho+y_{\rho}^{q_{C,\rho}};
\alpha_\rho+v_{\rho}^{i,1},...,\alpha_\rho+v_{\rho}^{i,
q_{\rho}^i})}\over
{\prod_{\rho\not=\rho_0}\prod_{i=1}^{q_{C,\rho}}(\alpha_\rho+y_{\rho}^i)}}[C]$$
in the quotient field $K$ of  $H^*(BG;{\Bbb Z}_2)={\Bbb
Z}_2[a_1,...,a_k]$.

Kosniowski and Stong \cite{ks2} gave the following formula.

\begin{thm} [Kosniowski-Stong] If $f(\alpha; x; x^i)$
is of degree less than or equal to $n$, then
$$f(\alpha; x; x^i)[M]=\sum_C{{f(\alpha; z, \alpha+y;
\alpha+v)}\over{\prod(\alpha+y)}}[C]$$ in $K$.
\end{thm}

Our interest in this paper is on the $({\Bbb Z}_2)^k$-actions
$(\Phi, M^n)$ with $w(C)=1$ for all connected components $C$ of
$F$. In this case, we may choose $z=0$ in Theorem 2.1 from the
splitting principle. Thus we have

\begin{cor}
Suppose that $(\Phi, M^n)$ is a $({\Bbb Z}_2)^k$-action on a
smooth closed manifold with $w(C)=1$ for all fixed connected
components $C$. Then for all $f(\alpha; x)$ of degree less than or
equal to $n$
$$f(\alpha; x)[M]=\sum_C{{f(\alpha; 0, \alpha+y)}\over{\prod(\alpha+y)}}[C]$$ in $K$.
\end{cor}

 Finally, we give the definition of the linear
independence property for the $p$-dimensional part $F^p$ of the
fixed point set of a $({\Bbb Z}_2)^k$-action $(\Phi, M^n)$.
Generally speaking, all normal dimensional sequences of components
of $F^p$ may not be distinct if $F^p$ is disconnected. However,
$(\Phi, M^n)$ must be cobordant to a $({\Bbb Z}_2)^k$-action such
that  all elements of the normal dimensional sequence set of its
$p$-dimensional part ${F'}^p$ are distinct and ${F'}^0$ is
possibly empty if $p=0$. In fact,  one may form a
 connected sum for those components in $F^p$ with the same normal
 dimensional sequence when $p>0$, and one may cancel pairs of
 components with the same normal dimensional sequence when $p=0$.
 This doesn't change the $({\Bbb Z}_2)^k$-action $(\Phi, M^n)$ up
 to equivariant cobordism. We say that $F^p$ possesses the
{\em linear independence property} if the following conditions are
satisfied:

1) When  $p>0$,   the  normal dimensional sequence set
$$\{\{q_{C, \rho}\}\vert C\text{ is a connected
component of ${F'}^p$}\}$$ of ${F'}^p$ possesses the property
$(*)$: all monomials
${1\over{\prod_{\rho\not=\rho_0}\alpha_{\rho}^{q_{C, \rho}}}}$ are
linearly independent in the quotient field of $H^*(BG;{\Bbb
Z}_2)={\Bbb Z}_2[a_1,...,a_k]$.

2) When $p=0$, ${F'}^0$ is either empty or nonempty, and for the
nonempty case, the normal dimensional sequence set of  ${F'}^0$
satisfies the property $(*)$.

Note: The definition for the linear independence property given
here seems to be weaker than that stated in \cite{l2} or
\cite{l3}. However, up to equivariant cobordism for any $({\Bbb
Z}_2)^k$-action without bundles of covering action, essentially
there is no genuine difference.

\section{The proof of Theorem 1.1}

 Suppose that $(\Phi, M^n)$ is a
smooth $({\Bbb Z}_2)^k$-action on a
 smooth closed manifold such that each $p$-dimensional part $F^p$ of  the fixed point set $F$
 possesses the linear independence property,  and $w(F)=1$ and
 $\dim M>2^k\dim F$.  Since the equivariant cobordism class of any
 $({\Bbb Z}_2)^k$-action is determined by its fixed data (see
 \cite{s}), it suffices to show that each connected component of $F$
 with its normal bundle bounds.

 First, let us consider all components of dimension greater than
zero.  Given a $p$-dimensional part $F^p=\sqcup_{j=1}^\ell C_j$
with $p>0$ and
 each $C_j$ connected, and let $$\sqcup_{j=1}^\ell
 \bigoplus_{\rho\not=\rho_0}\nu_{C_j, \rho}\longrightarrow C_j$$
 be the normal bundle to $F^p$ in $M$, and
 the normal dimensional sequence of each $C_j$ is $$\{q_{C_j,
 \rho}\vert \rho(\not=\rho_0)\in \text{Hom}(G,{\Bbb Z}_2)\}.$$
 Without loss of generality, one assumes that all normal dimensional
 sequences of $F^p$ are distinct.
 Now suppose inductively that each connected component of each
 $h$-dimensional part $F^h$ with its normal bundle bounds if
 $h>p$. Given partitions
 $$\omega_\rho=(i_1^\rho,...,i_{r_\rho}^\rho), \ \
 \rho\not=\rho_0$$
 with $\sum_{\rho\not=\rho_0}\vert\omega_\rho\vert=p$, choose
 $$f(\alpha;x)=\prod_{\rho\not=\rho_0}\{\sum
 [\prod_{\rho}(\alpha_\rho+x_1)]^{i_1^\rho}\cdots
 [\prod_{\rho}(\alpha_\rho+x_{r_\rho})]^{i_{r_\rho}^\rho}\}.$$
 Since $w(F)=1$, by Corollary 2.2 one has that

 \begin{eqnarray*}
 f(\alpha; 0, \alpha+y)&=&
 \prod_{\rho\not=\rho_0}s_{\omega_\rho}(y)\cdot
 (\prod_{\rho\not=\rho_0}\alpha_\rho)^s
  +\text{ terms of higher degree in $y$'s}
 \end{eqnarray*}
 where $s=\sum_{\rho\not=\rho_0}\vert \omega_\rho\vert=p$. Then,
 by induction and since deg$f(\alpha; 0, \alpha+y)$ is more than or equal to $p$ in
 $y$'s, $\sum_{h}{{f(\alpha, 0,
 \alpha+y)}\over{\prod(\alpha+y)}}[F^h]$ has no contribution for the component
 $F^h$ with $h\not=p$, i.e.,
 $$\sum_{h\not= p}{{f(\alpha; 0,
 \alpha+y)}\over{\prod(\alpha+y)}}[F^h]=0.$$
On the other hand, one has that the degree of $f(\alpha, x)$ in
$x$'s is
$$\deg f(\alpha; x)=\sum_{\rho\not=\rho_0}\vert \omega_\rho\vert +(2^k-1)\sum_{\rho\not=\rho_0}\vert
\omega_\rho\vert=2^k\sum_{\rho\not=\rho_0}\vert
\omega_\rho\vert=2^kp\leq 2^k\dim F<n$$ so by Corollary

\begin{eqnarray*}
0&=& f(\alpha; x)[M]\\
&=&\sum_{j=1}^\ell
{{\prod_{\rho\not=\rho_0}s_{\omega_\rho}(y)\cdot
 (\prod_{\rho\not=\rho_0}\alpha_\rho)^s +\text{ terms of higher degree in
 $y$'s}}\over{\prod_{\rho\not=\rho_0}\prod_{i=1}^{q_{C_j,
 \rho}}(\alpha_\rho+y^i_\rho)}}[C_j]\\
 &=&(\prod_{\rho\not=\rho_0}\alpha_\rho)^s\sum_{j=1}^\ell
 {{\prod_{\rho\not=\rho_0}s_{\omega_\rho}(y)}\over{\prod_{\rho\not=\rho_0}\alpha_{\rho}^{q_{C_j,
 \rho}}}}[C_j]
 \end{eqnarray*}
and thus
$$\sum_{j=1}^\ell
 {{\prod_{\rho\not=\rho_0}s_{\omega_\rho}(y)[C_j]}\over{\prod_{\rho\not=\rho_0}\alpha_{\rho}^{q_{C_j,
 \rho}}}}=0.$$
Since the $F^p$ possesses the linear independence property, one
obtains that for each $j$,
$$\prod_{\rho\not=\rho_0}s_{\omega_\rho}(y)[C_j]=0,$$
which means that $\bigoplus_{\rho\not=\rho_0}\nu_{C_j,
\rho}\longrightarrow C_j$ bounds. This completes the induction,
and shows that all components of dimension greater than zero with
their normal bundles bound. As for the 0-dimensional part $F^0$ of
the fixed point set, it is easy to see that either $F^0$ is empty
or the normal dimensional sequences of all isolated points in
$F^0$ must appear in pairs if $F^0$ is nonempty. This means that
$F^0$ with its normal bundle cobords away. Thus, $(\Phi, M^n)$
bounds equivariantly. This completes the proof of Theorem 1.1.
$\Box$

Now let us explain the Remarks (1) and (2) in Section 1 by the
following two examples.

 {\bf Example 1.} Begin with the involution $(T, {\Bbb R}P^2)$
given by $[x_0,x_1,x_2]\longmapsto [-x_0,x_1,x_2]$, which fixes
the disjoint union of a point and a real projective 1-space ${\Bbb
R}P^1$. Then the product
$$(T\underbrace{\times \cdots\times}_\ell T, {\Bbb
R}P^2\underbrace{\times\cdots\times}_\ell{\Bbb R}P^2)$$ of $\ell$
copies of $(T, {\Bbb R}P^2)$ forms a new involution, and its fixed
point set is $\bigsqcup_{i=0}^\ell {{\ell}\choose i}{\Bbb
R}P^1\underbrace{\times\cdots\times}_i{\Bbb R}P^1$, where ${\Bbb
R}P^1\underbrace{\times\cdots\times}_i{\Bbb R}P^1$ means a point
if $i=0$. This new involution is cobordant to an involution
$(\Phi_1, M_1^{2\ell})$ having fixed set $F=F^{\ell}\bigsqcup
F^{\ell-1}\bigsqcup\cdots\bigsqcup F^0$ with $\dim M^{2\ell}=2\dim
F$ and $w(F)=1$, where
$$F^p=
\begin{cases}
{\Bbb R}P^1\underbrace{\times\cdots\times}_p{\Bbb R}P^1 & \text{
if } {{\ell}\choose p}\not\equiv 0\mod 2\\
\text{empty} & \text{ if } {{\ell}\choose p}\equiv 0\mod 2.
\end{cases}
$$
Consider $M_1^{2\ell}\times M_1^{2\ell}$ with two involutions
$t_1=\text{twist}$ and $t_2=\Phi_1\times \Phi_1$. The fixed point
set of this $({\Bbb Z}_2)^2$-action $(\Phi_2, M_2^{4\ell})$ is the
fixed point set of $\Phi_1$ in the diagonal copy of $M_1^{2\ell}$
which is $F=F^{\ell}\bigsqcup F^{\ell-1}\bigsqcup\cdots\bigsqcup
F^0$, which has $w(F)=1$ and $\dim M_2^{4\ell}=2^2\dim F$.
Squaring this example gives examples for all $({\Bbb
Z}_2)^k$-actions. Actually,  if $(\Phi_{k-1},
M_{k-1}^{2^{k-1}\ell})$ is a $({\Bbb Z}_2)^{k-1}$-action fixing
$F=F^{\ell}\bigsqcup F^{\ell-1}\bigsqcup\cdots\bigsqcup F^0$, then
the twist  and the diagonal $({\Bbb Z}_2)^{k-1}$-action induced by
$\Phi_{k-1}$ on $M_{k-1}^{2^{k-1}\ell}\times
M_{k-1}^{2^{k-1}\ell}$ produce a $({\Bbb Z}_2)^{k}$-action
$(\Phi_{k}, M_{k}^{2^{k}\ell})$ whose fixed set is still
$F=F^{\ell}\bigsqcup F^{\ell-1}\bigsqcup\cdots\bigsqcup F^0$, and
$\dim M_k^{2^k\ell}=2^k\dim F$. Also, the linear independence for
the fixed point set is obvious since $F^p$ is connected for each
$p$. However, $(\Phi_{k}, M_{k}^{2^{k}\ell})$ is nonbounding for
every value of $\dim F=\ell$ and every $k$.

{\bf Example 2.} Consider the standard $({\Bbb Z}_2)^2$-action
$(\Phi_0, {\Bbb R}P^2)$ given by
$$[x_0, x_1, x_2]\longmapsto [x_0, g_1x_1, g_2x_2],$$ which fixes
three isolated points, where $(g_1,g_2)\in ({\Bbb Z}_2)^2$. Then
the diagonal action on  the product of $2\ell'$ copies of
$(\Phi_0, {\Bbb R}P^2)$ is also  a $({\Bbb Z}_2)^{2}$-action
denoted by $(\Phi, M^{4\ell'})$, and the fixed point set of this
action  is formed by $3^{2\ell'}$ isolated points.  Furthermore,
by using the construction as in Example 1 to $(\Phi, M^{2\ell'})$,
one may obtain a $({\Bbb Z}_2)^k$-action $(\Psi, M^{2^k\ell'})$,
which fixes $3^{2\ell'}$ isolated points. Now, the diagonal action
on the product of $(\Psi, M^{2^k\ell'})$ and  $(\Phi_k,
M_k^{2^k\ell})$ in Example 1 produces a $({\Bbb Z}_2)^{k}$-action
$(\Phi', M^{2^{k}(\ell+\ell')})$ fixing the disjoint union $F'$ of
$3^{2\ell'}$ copies of $F=F^{\ell}\bigsqcup
F^{\ell-1}\bigsqcup\cdots\bigsqcup F^0$. However, $(\Phi',
M^{2^{k}(\ell+\ell')})$ never bounds although $\dim
M^{2^{k}(\ell+\ell')}=2^k(\ell+\ell')> 2^{k}\dim F'$ for $\ell'>0$
and $w(F')=1$. This is because each $p$-dimensional part of $F'$
doesn't satisfy the linear independence property.

\section{ The case $\dim M=2^k\dim F$}

Suppose that  $(\Phi, M^n)$ is a $({\Bbb Z}_2)^k$-action with
 $w(F)=1$.
 Now let us
consider the case in which $\dim M^n=2^k\dim F.$

When $k=1$, one has

\begin{prop}
Let $(\Phi, M^n)$ is an involution with $w(F)=1$. If $\dim
M^n=2\dim F$, then $(\Phi, M^n)$ either bounds or is cobordant to
$(T\underbrace{\times\cdots\times}_{\dim F}T, {\Bbb
R}P^2\underbrace{\times\cdots\times}_{\dim F}{\Bbb R}P^2)$, where
the involution $(T, {\Bbb R}P^2)$ is given by
$[x_0,x_1,x_2]\longmapsto[-x_0,x_1,x_2]$.
\end{prop}

\begin{pf}
It suffices to show that if $(\Phi, M^n)$ is nonbounding, then it
is cobordant to $$(T\underbrace{\times\cdots\times}_{\dim F}T,
{\Bbb R}P^2\underbrace{\times\cdots\times}_{\dim F}{\Bbb R}P^2).$$
Suppose that $(\Phi, M^n)$ is nonbounding. Then the
$n/2$-dimensional component $F^{n/2}$ with normal bundle
$\nu^{n/2}$ must not bound.  For any a nondyadic partition
$\omega=(j_1,...,j_t)$ of $n$, one must have $t\leq n/2$ since
each $j_\alpha$ is not of the form $2^s-1$.  Take
$f_\omega(x)=\sum x_1(x_1+1)^{j_1-1}\cdots x_t(x_t+1)^{j_t-1}$,
one then has from \cite{ks1} that
$$s_{(j_1,...,j_t)}[M^n]=
\begin{cases}
0 & \text{ if } t\not=n/2\\
s_{(j_1-1,...,j_{n/2}-1)}(\nu^{n/2})[F^{n/2}] & \text{ if } t=n/2.
\end{cases} $$
However, $t=n/2$ forces $\omega$ to be $(2,\underbrace{...}_{n/2},
2)$, and so if $(\Phi, M^n)$ is nonbounding, then all
Stiefel-Whitney numbers of $M^n$ are zero except
$s_{(2,\underbrace{...}_{n/2}, 2)}[M^n]$. Thus one must have that
$M^n$ is cobordant to ${\Bbb
R}P^2\underbrace{\times\cdots\times}_{n/2}{\Bbb R}P^2$. By
\cite{ks1}, since the cobordism class of $M^n$ determines that of
$(\Phi, M^n)$, $(\Phi, M^n)$ is cobordant to
$$(T\underbrace{\times\cdots\times}_{\dim F}T,
{\Bbb R}P^2\underbrace{\times\cdots\times}_{\dim F}{\Bbb R}P^2).$$
\end{pf}

Proposition 4.1 shows that there is only one nonbounding
equivariant cobordism class of  involutions $(\Phi, M^n)$ with
$\dim M^n=2\dim F$ and $w(F)=1$.

For the case $k>1$,  generally there is no such  uniqueness result
 for nonbounding $({\Bbb Z}_2)^k$-actions with $w(F)=1$ and $\dim
M=2^k\dim F$ even if each $p$-dimensional part of the fixed point
set has the linear independence property. (Note: Pergher showed in
\cite{p} that if $(\Phi, M^n)$ is a $({\Bbb Z}_2)^k$-action with
$\dim M^n=2^k\dim F$ and $F$ being connected (here $F$ doesn't
necessarily  satisfy $w(F)=1$), then $(\Phi, M^n)$ is cobordant to
the twist on $F\underbrace{\times\cdots\times}_{2^k} F$. Thus, if
$F$ is connected with $w(F)=1$ and $\dim M=2^k\dim F$, then
$(\Phi, M^n)$ bounds equivariantly.) The question seems to be
quite complicated. Two kinds of examples are stated as follows.

First, let us look at the  $(\Phi_k, M_k^{2^k\ell})$ in Example 1.
By applying automorphisms of $({\Bbb Z}_2)^k$ to $\Phi_k$ to
switch normal representations around, one may obtain more $({\Bbb
Z}_2)^k$-actions, each of which is not cobordant to $(\Phi_k,
M_k^{2^k\ell})$, and  the disjoint union of any two such actions
also gives the example with the linear independence for the fixed
point set. Except for these examples, one may also construct other
examples as follows.

Consider the $({\Bbb Z}_2)^2$-action $(\Psi_0, {\Bbb R}P^4)$
defined by
$$t_1: [x_0,x_1,x_2,x_3,x_4]\longmapsto [-x_0,-x_1,x_2,x_3,x_4]$$
$$t_2: [x_0,x_1,x_2,x_3,x_4]\longmapsto [-x_0,x_1,-x_2,x_3,x_4],$$
whose fixed point set $F_0$ is the disjoint union of  three
isolated points $$p_1=[1,0,0,0,0],\ \ p_2=[0,1,0,0,0],\ \
p_3=[0,0,1,0,0]$$ and a real projective 1-space ${\Bbb R}P^1$.
Obviously, the 1-dimensional part ${\Bbb R}P^1$ of $F_0$ possesses
the linear independence property. Now let us show that the
0-dimensional part of $F_0$ also possesses the linear independence
property. Let $\rho_1, \rho_2, \rho_3$ be three irreducible
1-dimensional real $({\Bbb Z}_2)^2$-representations defined as
follows:
$$\rho_1(t_1)=-1, \rho_1(t_2)=1; \ \ \rho_2(t_1)=1,
\rho_2(t_2)=-1; \ \ \rho_3(t_1)=-1, \rho_3(t_2)=-1$$ where $t_1,
t_2$ are two generators of $({\Bbb Z}_2)^2$. Then it is easy to
see that the normal representations at three isolated points
$[1,0,0,0,0], [0,1,0,0,0], [0,0,1,0,0]$ are
$$\rho_1\rho_2\rho_3^2, \ \rho_1^2\rho_2\rho_3, \
\rho_1\rho_2^2\rho_3,$$ respectively. Furthermore, one knows that
the equivariant Euler classes at three isolated points $p_1, p_2,
p_3$ are
$$a_1a_2(a_1+a_2)^2,\ a_1^2a_2(a_1+a_2),\ a_1a_2^2(a_1+a_2),$$
respectively. Since
$${1\over{a_1a_2(a_1+a_2)^2}},\ {1\over{a_1^2a_2(a_1+a_2)}},\ {1\over{a_1a_2^2(a_1+a_2)}}$$
are linearly independent in the quotient field of ${\Bbb
Z}_2[a_1,a_2]$, one obtains that the 0-dimensional part of $F_0$
also possesses the linear independence property. Thus, $(\Psi_0,
{\Bbb R}P^4)$ is a $({\Bbb Z}_2)^2$-action with the linear
independence for the fixed point set. Furthermore, the diagonal
action on the product of $\ell$ copies of $(\Psi_0, {\Bbb R}P^4)$
is also a $({\Bbb Z}_2)^2$-action $(\Psi_2, N_2^{2^2\ell})$,  and
the fixed point set $F$ is the product of $\ell$ copies of $F_0$,
which has $w(F)=1$ and $\dim N_2^{2^2\ell}=2^2\dim F$. In
particular, the $r$-dimensional part $F^r$ of $F$ is a disjoint
union
$$\bigsqcup_{r_1+r_2+r_3+r=\ell}{{\ell}\choose{r_1,r_2,r_3,
r}}p_1^{r_1}p_2^{r_2}p_3^{r_3}({\Bbb R}P^1)^r$$ where the number
${{\ell}\choose{r_1,r_2,r_3, r}}={{\ell!}\over{r_1!r_2!r_3!r!}}$
is the multinomial coefficient, and $p_i^{r_i}$ means
$p_i\underbrace{\times\cdots\times}_{r_i}p_i$. Thus, $(\Psi_2,
N_2^{2^2\ell})$ is actually cobordant to a $({\Bbb Z}_2)^2$-action
such that the $r$-dimensional part of its fixed point set is a
union
\begin{eqnarray}
\bigsqcup_{\begin{Sb}r_1+r_2+r_3+r=\ell \\
{{\ell}\choose{r_1,r_2,r_3, r}}\equiv 1\mod 2 \end{Sb}}
p_1^{r_1}p_2^{r_2}p_3^{r_3}({\Bbb R}P^1)^r.
\end{eqnarray}
 With
this understood,  for a convenience in the following discussion
one will  regard the $r$-dimensional part $F^r$ of the fixed point
set of $(\Psi_2, N_2^{2^2\ell})$ as being the form (4.1). Also, it
is easy to see that all elements of the normal dimensional
sequence set of the form (4.1) are distinct and nonbounding.
Specifically, $({\Bbb R}P^1)^r$ with its nontrivial normal bundle
is nonbounding and the normal representation for
$p_1^{r_1}p_2^{r_2}p_3^{r_3}$ is
$\rho_1^{r_1+2r_2+r_3}\rho_2^{r_1+r_2+2r_3}\rho_3^{2r_1+r_2+r_3}=(\rho_1\rho_2\rho_3)^{\ell-r}\rho_1^{r_1}\rho_2^{r_2}
\rho_3^{r_3}$ and these are distinct representations.

Generally,
 $(\Psi_2, N_2^{2^2\ell})$  is not
 a $({\Bbb Z}_2)^2$-action with the linear independence for the
fixed point set. For example, taking $\ell=3$, one has that
 the 0-dimensional part $F^0$ of the fixed point set consists of nine isolated points
$$p_1^3, p_2^3, p_3^3, p_1^2p_2, p_1^2p_3, p_1p_2^2, p_2^2p_3,
p_1p_3^2, p_2p_3^2.$$ It is easy to see that the equivariant Euler
classes of these nine isolated points are $$a_1^3a_2^3(a_1+a_2)^6,
a_1^6a_2^3(a_1+a_2)^3, a_1^3a_2^6(a_1+a_2)^3,
a_1^4a_2^3(a_1+a_2)^5, a_1^3a_2^4(a_1+a_2)^5,$$
$$a_1^5a_2^3(a_1+a_2)^4, a_1^5a_2^4(a_1+a_2)^3,
a_1^3a_2^5(a_1+a_2)^4, a_1^4a_2^5(a_1+a_2)^3,$$ respectively. Now
let us look at four isolated points $p_1^2p_2, p_1^2p_3,
p_1p_3^2,p_2p_3^2$. One has that
\begin{eqnarray*}
&&{1\over{a_1^4a_2^3(a_1+a_2)^5}}+{1\over{a_1^3a_2^4(a_1+a_2)^5}}+{1\over{a_1^3a_2^5(a_1+a_2)^4}}+
{1\over{a_1^4a_2^5(a_1+a_2)^3}}\\
&=&{{a_1^2a_2^3(a_1+a_2)+a_1^3a_2^2(a_1+a_2)+
a_1^3a_2(a_1+a_2)^2+a_1^2a_2(a_1+a_2)^3}\over{a_1^6a_2^6(a_1+a_2)^6}}\\
&=&0
\end{eqnarray*}
in the quotient field of ${\Bbb Z}_2[a_1,a_2]$, so linear
independence fails.

Now let us discuss when $(\Psi_2, N_2^{2^2\ell})$ is  a $({\Bbb
Z}_2)^2$-action with the linear independence for the fixed point
set. Note that ${{\ell}\choose{r_1,r_2,r_3,r}}= {{\ell}\choose
r}{{\ell-r}\choose{r_1,r_2,r_3}}$, so
${{\ell}\choose{r_1,r_2,r_3,r}}\equiv 1 \mod 2$ if and only if $
{{\ell}\choose r}\equiv 1\mod 2$ and $
{{\ell-r}\choose{r_1,r_2,r_3}}\equiv 1 \mod 2$.

For the $r$-dimensional part $F^r$ with ${{\ell}\choose r}\equiv
1\mod 2$, since each component $p_1^{r_1}p_2^{r_2}p_3^{r_3}({\Bbb
R}P^1)^r$ with ${{\ell-r}\choose{r_1,r_2,r_3}}\equiv 1\mod 2$
contains $({\Bbb R}P^1)^r$ as a factor with the same normal
bundle, and since the equivariant Euler class of
$p_1^{r_1}p_2^{r_2}p_3^{r_3}({\Bbb
R}P^1)^r\vert_{p_1^{r_1}p_2^{r_2}p_3^{r_3}}$ restricted to
$p_1^{r_1}p_2^{r_2}p_3^{r_3}$ is
$$a_1^{r_1+2r_2+r_3}a_2^{r_1+r_2+2r_3}(a_1+a_2)^{2r_1+r_2+r_3}, $$
the linear independence of $F^r$ is equivalent to  that of
$$B_r=\{{1\over
{a_1^{r_1+2r_2+r_3}a_2^{r_1+r_2+2r_3}(a_1+a_2)^{2r_1+r_2+r_3}}}\vert
r_1+r_2+r_3=\ell-r, {{\ell-r}\choose{r_1,r_2,r_3}}\equiv 1\mod
2\}$$ in the quotient field of ${\Bbb Z}_2[a_1,a_2]$. Thus one has

{\bf Claim 1.} {\em The linear independence of $F^r$ is equivalent
to that of $B_r$ in the quotient field of ${\Bbb Z}_2[a_1,a_2]$.}

On the other hand, all elements of $B_r$ are actually all nonzero
monomials formed by ${1\over {a_1a_2(a_1+a_2)^2}}$, ${1\over
{a_1^2a_2(a_1+a_2)}}$, ${1\over {a_1a_2^2(a_1+a_2)}}$ of the
expression  of $$({1\over {a_1a_2(a_1+a_2)^2}}+{1\over
{a_1^2a_2(a_1+a_2)}}+{1\over {a_1a_2^2(a_1+a_2)}})^{\ell-r}$$ over
${\Bbb Z}_2$.  Since
$$({1\over
{a_1a_2(a_1+a_2)^2}}+{1\over {a_1^2a_2(a_1+a_2)}}+{1\over
{a_1a_2^2(a_1+a_2)}})^{\ell-r}={{(a_1a_2+a_2(a_1+a_2)+a_1(a_1+a_2))^{\ell-r}}\over
{(a_1a_2(a_1+a_2))^{2\ell-2r}}},$$ the problem is reduced to
determining when all elements of
$$D_r^{\ell-r}=\{a_1^{r_1+r_3}a_2^{r_1+r_2}(a_1+a_2)^{r_2+r_3}\vert
r_1+r_2+r_3=\ell-r, {{\ell-r}\choose{r_1,r_2,r_3}}\equiv 1\mod
2\}$$ are linearly independent in ${\Bbb Z}_2[a_1,a_2]$, where
$D_r^{\ell-r}$ just consists of all nonzero monomials formed by
$a_1a_2, a_1(a_1+a_2), a_2(a_1+a_2)$ of the expression  of
$(a_1a_2+a_2(a_1+a_2)+a_1(a_1+a_2))^{\ell-r}$ over ${\Bbb Z}_2$.

{\bf Definition.} Let $\ell=2^{s_1}+\cdots+2^{s_u}$ with
$s_1<\cdots<s_u$  be the 2-adic expansion of $\ell$. Then  $\ell$
has the {\em  gap property} if $u=1$ or $u>1$ and $s_{i+1}-s_i>1$
for each $1\leq i\leq u-1$.

{\bf Claim 2.} {\em  $\ell$ has the  property that all elements of
$D_r^{\ell-r}$ are linearly independent in ${\Bbb Z}_2[a_1,a_2]$
for each $0\leq r\leq \ell$ with ${{\ell}\choose r}\equiv 1\mod 2$
if and only if $\ell$ has the gap property.}
\begin{pf}
 Suppose $\ell$ has the  gap property. One  uses induction on
 $u$. If $u=1$, then  $\ell$ is of the form $2^s$ so
$${{\ell}\choose{r_1, r_2, r_3, r}}\equiv 0 \mod 2$$
except that one of $r_1, r_2, r_3, r$ is equal to $\ell$. This
means that $r=0$ or $\ell$. When $r=0$, $D_0^\ell$ contains three
nonzero terms $a_1^{\ell}a_2^{\ell}, a_1^{\ell}(a_1+a_2)^{\ell},
a_2^{\ell}(a_1+a_2)^{\ell}$, which are obviously linearly
independent in ${\Bbb Z}_2[a_1,a_2]$. When $r=\ell$, $D_\ell^0$
contains only the term 1, which is also linearly independent in
${\Bbb Z}_2[a_1,a_2]$. If $u\leq v$, suppose inductively that for
each $0\leq r\leq \ell$ with ${{\ell}\choose r}\equiv 1\mod 2$,
all elements of $D_r^{\ell-r}$ are linearly independent in ${\Bbb
Z}_2[a_1,a_2]$. Consider the case $u=v+1$, for  $r>0$ with
${{\ell}\choose r}\equiv 1\mod 2$, since ${{\ell}\choose r}\equiv
1\mod 2$, one knows that the 2-adic expansion of $r$ is  part of
that of $\ell$, so $\ell-r$ has the  gap property and the number
of  terms of the 2-adic expansion of $\ell-r$ is at most $v$.
Thus, by induction, if $r>0$ with ${{\ell}\choose r}\equiv 1\mod
2$, then all elements of $D_r^{\ell-r}$ are linearly independent
in ${\Bbb Z}_2[a_1,a_2]$. As for $r=0$, write
$\ell=\ell'+2^{s_{v+1}}$, one then has that
\begin{eqnarray*}
&&\{a_1a_2+a_1(a_1+a_2)+a_2(a_1+a_2)\}^\ell\\
&=&\{a_1a_2+a_1(a_1+a_2)+a_2(a_1+a_2)\}^{\ell'}\\
&&\times(a_1^{2^{s_{v+1}}}
a_2^{2^{s_{v+1}}}+a_1^{2^{s_{v+1}}}(a_1+a_2)^{2^{s_{v+1}}}+a_2^{2^{s_{v+1}}}(a_1+a_2)^{2^{s_{v+1}}})
\end{eqnarray*}
so
$$D_0^\ell=\{a_1^{2^{s_{v+1}}}
a_2^{2^{s_{v+1}}}\}\times D_0^{\ell'}\bigsqcup \{a_1^{2^{s_{v+1}}}
(a_1+a_2)^{2^{s_{v+1}}}\}\times
D_0^{\ell'}\bigsqcup\{a_2^{2^{s_{v+1}}}
(a_1+a_2)^{2^{s_{v+1}}}\}\times D_0^{\ell'}.$$
Let
\begin{eqnarray*}
0&=&a_1^{2^{s_{v+1}}} a_2^{2^{s_{v+1}}}\sum_{X\in
D_0^{\ell'}}l_X^{(1)}X+a_1^{2^{s_{v+1}}}
(a_1+a_2)^{2^{s_{v+1}}}\sum_{X\in
D_0^{\ell'}}l_X^{(2)}X+a_2^{2^{s_{v+1}}}
(a_1+a_2)^{2^{s_{v+1}}}\sum_{X\in D_0^{\ell'}}l_X^{(3)}X\\
&=& a_1^{2^{s_{v+1}+1}}\sum_{X\in D_0^{\ell'}}l_X^{(2)}X+
a_2^{2^{s_{v+1}+1}}\sum_{X\in
D_0^{\ell'}}l_X^{(3)}X+a_1^{2^{s_{v+1}}}
a_2^{2^{s_{v+1}}}\sum_{X\in
D_0^{\ell'}}(l_X^{(1)}+l_X^{(2)}+l_X^{(3)})X
\end{eqnarray*}
where $l_X^{(1)}, l_X^{(2)}, l_X^{(3)}\in {\Bbb Z}_2$. Since
$s_{v+1}-s_i>1$ for any $i\leq v$, this forces
$$\sum_{X\in D_0^{\ell'}}l_X^{(1)}X=0, \sum_{X\in D_0^{\ell'}}l_X^{(2)}X=0, \sum_{X\in
D_0^{\ell'}}l_X^{(3)}X=0.$$ By induction, one has that all
$l_X^{(1)},l_X^{(2)},l_X^{(3)}$ are zero, so all elements of
$D_0^\ell$ are linearly independent in ${\Bbb Z}_2[a_1,a_2]$. This
completes the induction.

Now suppose  $\ell$ does not have the gap  property. There is at
least one $i$ with $s_{i+1}-s_i=1$. Take
$r_0=\ell-(2^{s}+2^{s+1})$ where $s_i=s$ and  then
$D_{r_0}^{\ell-r_0}$ contains  the nine terms  of
$(a_1a_2+a_1(a_1+a_2)+a_2(a_1+a_2))^{3\cdot 2^s}$ or the $2^s$
powers of the elements of $D_0^3$. Then
$$0=(a_1^2a_2^3(a_1+a_2)+a_1^3a_2^2(a_1+a_2)+a_1^3a_2(a_1+a_2)^2+a_1^2a_2(a_1+a_2)^3)^{2^s}$$
gives a linear dependence relation, exactly as in the case
$\ell=3$.
\end{pf}

Combining Claims 1 and 2, one has

{\bf Fact.} {\em  $(\Psi_2, N_2^{2^2\ell})$ is  a $({\Bbb
Z}_2)^2$-action with the linear independence for the fixed point
set if and only if  $\ell$ has the  gap property.}

 Now suppose that $\ell$ has the  gap property. Using the construction in Example 1, one can
obtain a $({\Bbb Z}_2)^k$-action $(\Psi_k, N_k^{2^k\ell})$, and
the fixed point set $F$ is still the same as that of $(\Psi_2,
N_2^{2^2\ell})$. Obviously,
 $(\Psi_k, N_k^{2^k\ell})$ is
nonbounding  and is   a $({\Bbb Z}_2)^k$-action  with $w(F)=1$ and
$\dim N_k^{2^k\ell}=2^k\dim F$. The linear independence of each
$r$-dimensional part $F^r$ of the fixed point set of $(\Psi_k,
N_k^{2^k\ell})$ follows from the following lemma 4.1, so $(\Psi_k,
N_k^{2^k\ell})$ is also a $({\Bbb Z}_2)^k$-action with the linear
independence for the fixed point set.

Both $(\Phi_k, M_k^{2^k\ell})$ and $(\Psi_k, N_k^{2^k\ell})$ give
two different kinds of examples since $M_k^{2^k\ell}$ is never
cobordant to $N_k^{2^k\ell}$.

\begin{lem}
Suppose that $(\Phi, M^n)$ is a $({\Bbb Z}_2)^k$-action satisfying
linear independence for the fixed point set. Then $(\Psi,
M^n\times M^n)$, which is the $({\Bbb Z}_2)^{k+1}$-action formed
by the diagonal action $\Phi\times\Phi$ and twist on $M^n\times
M^n$, also satisfies the linear independence for its fixed point
set.
\end{lem}
\begin{pf}
The fixed point set of the $({\Bbb Z}_2)^{k+1}$-action $\Psi$ on
$M^n\times M^n$ is the same as the fixed point set of $(\Phi,
M^n)$, but the normal representation changes. Along a fixed
connected component $C^p$ the tangent bundle of $M^n$ restricts to
a $({\Bbb Z}_2)^k$-representation
$$\prod_{\rho\in\text{Hom}(({\Bbb Z}_2)^k,{\Bbb
Z}_2)}\rho^{q_{C,\rho}}.$$ Note that $\rho_0^{q_{C,\rho_0}}$
corresponds to the tangent bundle of $C^p$ so $q_{C,\rho_0}=p$.
Let $\Lambda$ and $\bar{\Lambda}$ be two subsets of
$\text{Hom}(({\Bbb Z}_2)^{k+1},{\Bbb Z}_2)$ such that both
$\Lambda$ and $\bar{\Lambda}$ are isomorphic to $\text{Hom}(({\Bbb
Z}_2)^{k},{\Bbb Z}_2)$ as ${\Bbb Z}_2$ vector spaces, and each
$\delta_\rho$ in $\Lambda$ is $\rho\in\text{Hom}(({\Bbb
Z}_2)^k,{\Bbb Z}_2)$ on $({\Bbb Z}_2)^k$ and 1 on the new ${\Bbb
Z}_2$ generator $t_{k+1}$ and each $\bar{\delta}_\rho$ in
$\bar{\Lambda}$ is $\rho\in\text{Hom}(({\Bbb Z}_2)^k,{\Bbb Z}_2)$
on $({\Bbb Z}_2)^k$ and $-1$ on the new ${\Bbb Z}_2$ generator
$t_{k+1}$.  Then, the tangent bundle of $M^n\times M^n$ along
$C^p$ restricts to a $({\Bbb Z}_2)^{k+1}$-representation
$$\prod_{\delta_\rho\in\Lambda}\delta_{\rho}^{q_{C,\rho}}\prod_{\bar{\delta}_\rho\in
\bar{\Lambda}}\bar{\delta}_{\rho}^{q_{C,\rho}}.$$ In particular,
the normal $({\Bbb Z}_2)^{k+1}$-representation of $C^p$ in
$M^n\times M^n$ is
$$\prod_{\delta_\rho\in\Lambda, \delta_\rho\not=\delta_{\rho_0}}\delta_{\rho}^{q_{C,\rho}}\prod_{\bar{\delta}_\rho\in
\bar{\Lambda}}\bar{\delta}_{\rho}^{q_{C,\rho}}.$$ As ${\Bbb Z}_2$
vector spaces, $\text{Hom}(({\Bbb Z}_2)^{k+1},{\Bbb Z}_2)\cong
H^1(B({\Bbb Z}_2)^{k+1};{\Bbb Z}_2)$, and  $\Lambda\cong
\bar{\Lambda}\cong \text{Hom}(({\Bbb Z}_2)^k,{\Bbb Z}_2)$. Without
loss of generality, one may assume that each nontrivial element
$\delta_\rho$ in $\Lambda$ corresponds to a polynomial of degree
one in $a_1,...,a_k$, denoted by $\alpha_\rho\in H^1(B({\Bbb
Z}_2)^{k+1};{\Bbb Z}_2)$,  and of course, the trivial element
$\delta_{\rho_0}$ in $\Lambda$ corresponds to zero element of
$H^1(B({\Bbb Z}_2)^{k+1};{\Bbb Z}_2)$, where $H^*(B({\Bbb
Z}_2)^{k+1};{\Bbb Z}_2)={\Bbb Z}_2[a_1,...,a_{k+1}]$. Then one
sees that each $\bar{\delta}_\rho$ in $\bar{\Lambda}$ corresponds
to the polynomial $\alpha_\rho+a_{k+1}$ in $H^1(B({\Bbb
Z}_2)^{k+1};{\Bbb Z}_2)$. Thus, the equivariant Euler class of the
normal $({\Bbb Z}_2)^{k+1}$-representation of $C^p$ is
$$\prod_{\delta_\rho\in \Lambda,
\delta_\rho\not=\delta_{\rho_0}}\alpha_{\rho}^{q_{C, \rho}}\cdot
a_{k+1}^p\cdot\prod_{\delta_\rho\in \Lambda,
\delta_\rho\not=\delta_{\rho_0}}(\alpha_{\rho}+a_{k+1})^{q_{C,
\rho}}=a_{k+1}^p\prod_{\delta_\rho\in \Lambda,
\delta_\rho\not=\delta_{\rho_0}}[\alpha_{\rho}(\alpha_\rho+a_{k+1})]^{q_{C,
\rho}}.$$ Now consider the $p$-dimensional part $F^p$ of the fixed
point set of $(\Psi, M^n\times M^n)$. Let
$$\sum_{C^p\subset F^p}{{l_C}\over{a_{k+1}^p\prod_{\delta_\rho\in \Lambda,
\delta_\rho\not=\delta_{\rho_0}}[\alpha_{\rho}(\alpha_\rho+a_{k+1})]^{q_{C,
\rho}}}}=0$$ where $l_C\in {\Bbb Z}_2$. Since
$$\sum_{C^p\subset F^p}{{l_C}\over{a_{k+1}^p\prod_{\delta_\rho\in \Lambda,
\delta_\rho\not=\delta_{\rho_0}}[\alpha_{\rho}(\alpha_\rho+a_{k+1})]^{q_{C,
\rho}}}}={1\over{a_{k+1}^p}}\sum_{C^p\subset
F^p}{{l_C}\over{\prod_{\delta_\rho\in \Lambda,
\delta_\rho\not=\delta_{\rho_0}}[\alpha_{\rho}(\alpha_\rho+a_{k+1})]^{q_{C,
\rho}}}}$$ one has $$\sum_{C^p\subset
F^p}{{l_C}\over{\prod_{\delta_\rho\in \Lambda,
\delta_\rho\not=\delta_{\rho_0}}[\alpha_{\rho}(\alpha_\rho+a_{k+1})]^{q_{C,
\rho}}}}=0.$$ Since $\prod_{\delta_\rho\in \Lambda,
\delta_\rho\not=\delta_{\rho_0}}\alpha_{\rho}^{q_{C,
\rho}}\not=0$, after reducing modulo $a_{k+1}$ one has that
$$\sum_{C^p\subset F^p}{{l_C}\over{\prod_{\delta_\rho\in \Lambda,
\delta_\rho\not=\delta_{\rho_0}}\alpha_{\rho}^{2q_{C,
\rho}}}}=(\sum_{C^p\subset F^p}{{l_C}\over{\prod_{\delta_\rho\in
\Lambda, \delta_\rho\not=\delta_{\rho_0}}\alpha_{\rho}^{q_{C,
\rho}}}})^2=0.$$ Squaring is a monomorphism in the quotient field
of ${\Bbb Z}_2[a_1,...,a_{k+1}]$, and thus one has
$$\sum_{C^p\subset F^p}{{l_C}\over{\prod_{\delta_\rho\in \Lambda,
\delta_\rho\not=\delta_{\rho_0}}\alpha_{\rho}^{q_{C, \rho}}}}=0$$
so $l_C=0$ for all $C^p$ in $F^p$ since $\prod_{\delta_\rho\in
\Lambda, \delta_\rho\not=\delta_{\rho_0}}\alpha_{\rho}^{q_{C,
\rho}}$ is actually identified with the equivariant Euler class of
the normal $({\Bbb Z}_2)^k$-representation of $C^p$ as a fixed
component of $(\Phi, M^n)$.
\end{pf}

}

\end{document}